\documentclass[12pt]{amsart}
\usepackage{geometry}
\usepackage{amsmath, amssymb, amsthm, amsfonts, amsxtra, mathrsfs, mathtools} 
\usepackage[intoc]{nomencl}
\makenomenclature
\usepackage{import}
\usepackage{color}
\usepackage[colorlinks=true, allcolors=blue]{hyperref}
\usepackage{stmaryrd}
\usepackage{tikz-cd}
\usepackage{tikz}
\usepackage{tikzsymbols}
\usepackage{enumitem}
\usetikzlibrary{matrix}
\usepackage[ruled,vlined]{algorithm2e}
\usepackage[mathscr]{euscript}

\input xy
\xyoption{all}

\numberwithin{equation}{section}

\newtheorem{thm}{Theorem}[section]
\newtheorem{lemma}{Lemma}[section]
\newtheorem{proposition}[thm]{Proposition}
\newtheorem{corollary}{Corollary}[section]

\theoremstyle{definition}
\newtheorem{conjecture}[thm]{Conjecture}
\newtheorem{definition}[thm]{Definition}
\newtheorem{solution}{Solution}
\newtheorem{example}[thm]{Example}

\newtheorem{question}[thm]{Question}

\newtheorem{fact}[thm]{Fact}

\def\G{\mathbf{G}}

\def\P{\mathbf{P}}

\def\calH{\mathcal{H}}

\def\calL{\mathcal{L}}

\renewcommand{\a}{\alpha}

\renewcommand{\l}{\lambda}

\newcommand{\s}{\sigma}

\newcommand\bl{\bullet}

\newcommand\eflag{E_{\bullet}}

\newcommand\Fl{\mathrm{Fl}}
\newcommand\fflag{F_{\bl}}
\newcommand\Fr{\textup{Fr}}

\newcommand{\GL}{\textup{GL}}
\newcommand{\Gr}{\textup{Gr}}

\newcommand\id{\textup{id}}

\newcommand\inv{\operatorname{inv}}

\newcommand\Pic{\textup{Pic}}

\newcommand\Spec{\operatorname{Spec}}

\newcommand{\pflag}{P_{\bullet}}

\newcommand\qflag{Q_{\bullet}}

\newcommand{\vflag}{V_{\bl}}
\newcommand{\vsa}{a_{\bl}} 
\newcommand{\vsb}{b_{\bl}} 

\newlength\mylen
\newlist{mycases}{enumerate}{1}
\setlist[mycases,1]{label=\textbf{Case~\arabic*.}, 
  labelwidth=\dimexpr-\mylen-\labelsep\relax,leftmargin=0pt,align=right}

\newcommand{\repeatable}[2]{%
  \label{#1}\global\@namedef{repeatable@#1}{#2}#2
}
\renewcommand{\repeat}[1]{%
  \@ifundefined{repeatable@#1}{NOT FOUND}{$\@nameuse{repeatable@#1}$}%
  ~\ref{#1}}
\makeatother

\title{Relative Bott-Samelson Varieties}
\author{Shiyue Li}\address{Shiyue Li, Department of Mathematics, Brown University, Providence, RI 02912}\email{shiyue\_li@brown.edu}
\date{\today}

\date{\today}

\begin{document}
\maketitle

\begin{abstract}
    We prove that, defined with respect to \textit{versal} flags, the product of two relative Bott-Samelson varieties over the flag bundle is a resolution of singularities of a \textit{relative Richardson variety}.
    This result generalizes Brion's resolution of singularities of Richardson varieties to the relative setting \cite{Brion2005}. It reflects the phenomenon that the local geometry of a relative Richardson variety is completely governed by the two intersecting relative Schubert varieties, studied in ~\cite{chan2019relative}.
    We also prove an analogous theorem in the case of relative Grassmannian Richardson varieties, thereby furnishing a resolution of singularities for the Brill-Noether variety with imposed ramification on twice-marked elliptic curves. 
\end{abstract}

\section{Introduction}
    A \textit{Bott-Samelson variety} is an iterated tower of $\P^1$-bundles or Grassmannian bundles.
    It is a subscheme of a product of flag varieties or Grassmannians, defined by  incidence relations prescribed by a reduced decomposition of a permutation $\s$ with respect to a fixed flag. A Bott-Samelson variety can serve as a resolution of singularities for a Schubert variety defined by $\s$ with respect to the same fixed flag. 
    One example of a Bott-Samelson resolution is the subvariety of the product of two Grassmannians $\G(0, 3) \times \G(1, 3)$ parametrizing pairs $(p, L)$ where $L$ intersects a fixed line in $\P^3$ and $p$ is contained in the intersection of $L$ with the fixed line. It admits a proper birational map to the Schubert variety in $\G(1, 3)$ that parametrizes lines intersecting a fixed line in $\P^3$; the birational map simply forgets the point $p$ in each pair $(p, L)$. 
    
    Furthermore, the intersection of two flag Schubert varieties defined with respect to transverse fixed flags is a \textit{Richardson variety} and is well-known to have rational singularities. In \cite{Brion2005}, Brion showed that, given a Richardson variety, the product of the Bott-Samelson resolutions for the two intersecting flag Schubert varieties over the complete flag varieties is a resolution of singularities for the Richardson variety.
    
    We are interested in generalizing Brion's result to a relative setting, motivated by Brill-Noether theory with imposed ramifications studied in \cite{osserman2011simple, Chan2017EulerCO, cop}. We ask the following question. 
    \begin{question}
    \label{qtn:main}
    Are there Bott-Samelson type resolutions of singularities for the intersection of two \emph{relative} Schubert varieties, under certain conditions on the defining flags? 
    \end{question}

    In recent work ~\cite{chan2019relative}, Chan and Pflueger showed that any local property of the intersection of two relative Schubert varieties is completely controlled by that of the relative Schubert varieties, under the additional condition that the defining flags are \textit{versal}. 
    Versality generalizes transversality: in the case when the base scheme is a point, versality recovers transversality. In the case when the base scheme is a $1$-dimensional scheme, versality implies that the flags are transverse in most fibers but become almost transverse (i.e.  exactly one pair of complementary dimensions of the two flags has a $1$-dimensional intersection and the remaining pairs intersect trivially; see  \cite{cop}) over finitely many reduced points. Their result generalizes \cite{KWY2013singularities}.
    
    In this work, we prove the product of two relative Bott-Samelson varieties over the flag bundle or the Grassmann bundle, with respect to versal defining flags, is a resolution of singularities for a relative Richardson variety in the flag bundle or the Grassmann bundle respectively, giving a positive answer to Question \ref{qtn:main}. 
    
    Throughout, let $k$ be an algebraically closed field of characteristic $0$. We state the main result as follows, postponing relevant definitions. 
    \begin{thm}
    \label{thm:main-2}
        Let $n > 0$. 
        Let $S$ be a smooth finite-type $k$-scheme that carries a rank-$n$ vector bundle $\calH$. 
        Let $\s, \tau \in S_n$ with lengths $\ell$ and $\ell'$ respectively. Let $\rho_{\s}$ and $\rho_{\tau}$ be reduced decompositions for $\s$ and $\tau$ respectively.
        Let $p, q: S \to \Fl(\calH)$  be versal sections of the associated flag bundle of $\calH$. Let $X_{\s}(p)$ and $X_{\tau}(q)$ be the relative Schubert varieties defined with respect to $p$ and $q$.
        Let $Z_{\rho_{\s}}(p)$ and $Z_{\rho_{\tau}}(q)$ be the relative Bott-Samelson varieties defined with respect to $p$ and $q$. 
        The fiber product
        \[Z_{\rho_{\s}}(p) \times_{\Fl(\calH)} Z_{\rho_{\tau}}(q)\]
        \begin{enumerate}
            \item is smooth; 
            \item and the projection to the relative Richardson variety \[
            Z_{\rho_{\s}}(p) \times_{\Fl(\calH)} Z_{\rho_{\tau}}(q) \to X_{\s}(p) \cap X_{\tau}(q)
            \] is proper and birational.
        \end{enumerate}
        Therefore, $Z_{\rho_{\s}}(p) \times_{\Fl(\calH)} Z_{\rho_{\tau}}(q)$ is a resolution of singularities of the relative Richardson variety.
    \end{thm}
    This reflects the nature that local geometric properties of relative Richardson varieties factor and generalizes Brion's result in the absolute case without explicitly invoking Kleiman's transversality theorem. 
    
    Our second result is an application of the above theorem in the case of relative Grassmannian Richardson varieties to Brill-Noether theory. 
    Let $d > r \ge 0$, and let $(E; P, Q)$ be a twice-marked elliptic curve where $P - Q$ is not a torsion point of order weakly less than $d$ in $\Pic^{0}(E)$. Denote by $\vsa$ and $\vsb$ sequences
    $0 \le a_1 < \cdots < a_r < a_{r+1}\le d$ and $0 \le b_1 < \cdots < b_r < b_{r+1} \le d$.
    The moduli space of linear series of projective rank $r$ with imposed ramification and of degree $d$ at $P, Q$ prescribed by $\vsa, \vsb$ respectively is the subscheme of the classical Brill-Noether variety $G^{r}_{d}(E)$
    \[
        G^{r}_{d}(E, (P, \vsa), (Q, \vsb)) \subseteq G^r_d(E), 
    \] and is the intersection of two relative Grassmannian Schubert varieties $G^r_d(E, (P, \vsa))$, $G^{r}_d(E, (Q, \vsb))$ defined with respect to versal flags; see \cite{cop, chan2019relative}.
    Let $\lambda = (\lambda_i)$ and $\lambda' = (\lambda'_{i'})$ be partitions where \begin{align*}\lambda_i &= a_{r + 1 - (i-1)} - (r+1 - i),\\\lambda'_{i'} &= b_{r + 1 - (i'-1)} - (r+1-i'),
    \end{align*} for $i, i' \in [r+1]$. Let $Z_{\lambda}(p)$ and $Z_{\lambda'}(q)$ be the relative Bott-Samelson resolutions for the relative Grassmannian Schubert varieties $\Gr^r_{d}(E, (P, \vsa))$ and $\Gr^r_{d}(E, (Q, \vsb))$ respectively; see Section \ref{sec:product-grass}. Then we have the following theorem. 
    \begin{thm}
    \label{thm:main-brill-noether-resing}
        The product of Bott-Samelson resolutions over the Brill-Noether variety $G^{r}_d(E)$
        \[
        Z_{\lambda}(p) \times_{G^{r}_{d}(E)} Z_{\lambda'}(q)
        \] is a resolution of singularities of $G^{r}_{d}(E, (P, \vsa), (Q, \vsb))$. 
    \end{thm}
    The resolution $Z_{\lambda}(p) \times_{G^{r}_{d}(E)} Z_{\lambda'}(q)$ has a natural modular interpretation: it parametrizes pairs of partial flags of linear series on $E$ of projective rank no greater than $r$, with ranks specified by $\vsa$ and $\vsb$ and with agreeing top-rank pieces. 
    At the end, we formulate two open conjectures regarding a generalization of such construction to partial defining flags and Bott-Samelson type resolutions of singularities for Brill-Noether varieties on twice-marked curves of higher genus. 

    \subsection*{Acknowledgements} The author would like to express thanks to Melody Chan for suggesting the problem and guidance throughout.  The author is also grateful to Dan Abromovich, Shamil Asgarli, Giovanni Inchiostro for helpful discussions.

\section{Background}

\subsection{Flag varieties, Schubert varieties and Bruhat stratification}
\label{sec:flag-schubert-varieties}
    We start with an introduction to flag varieties.
    The \textit{complete flag variety} of a fixed $n$-dimensional vector space $H$ is   
    \[
    \Fl(H) = \{\fflag = (\{0\} \subsetneq F_1 \subsetneq \cdots \subsetneq F_n = H): \dim F_i = i \text{ for all } i = 1, \ldots, n\}.
    \]
    Given a permutation $\s \in S_n$ and a flag $\fflag \in \Fl(H)$, the \textit{Schubert cell} $X^{\circ}_{\s, \fflag}$ is the subscheme in $\Fl(H)$ defined by incidence relations imposed by $\s$ with respect to the flag $\fflag$ and is set-theoretically 
    \[
    \left\{\vflag \in \Fl(H): \dim V_i \cap F_j = \#\{m: m \le i, \s(m) \le j\} \text{ for all } i, j \right\}.
    \] 
    The \textit{Schubert variety} $X_{\s, \fflag}$
    is the Zariski closure of the Schubert cell $X^{\circ}_{\s, \fflag}$ in $\Fl(H)$ and is set-theoretically
    \[
    \left\{\vflag \in \Fl(H): \dim V_{i} \cap F_j \ge \#\{m: m \le i, \s(m) \le j\} \text{ for all } i, j\right\}. 
    \] These incidence relations are indeed determinantal and thus bestow Schubert varieties with scheme structures. 
    Furthermore, if a flag $\fflag$ in $\Fl(H)$ is fixed, then the flag variety $\Fl(H)$ is the disjoint union of all Schubert cells defined with respect to $\fflag$, indexed by $S_n$. Formally,
    \[
    \Fl(H) = \coprod_{\s \in S_n} X^{\circ}_{\s, \fflag}.
    \]
    The set of Schubert varieties $\{X_{\s, \fflag}\}_{\s \in S_n}$ gives a stratification called the \textit{Bruhat stratification} on the flag variety $\Fl(H)$ and the containment relations of the strata yield a partial order called the \textit{Bruhat order} on these strata, as well as on $S_n$.
    It is always possible to write any permutation as a product of adjacent transpositions and if such a product has length 
    \[
    \inv(\s) = \#\{(i, j): \s(i) > \s(j) \text{ for all } 1 \le i < j \le n\},
    \] it is called the reduced decomposition of $\s$. The Bruhat order is the transitive closure of the relation:
    $\s \le \tau$ if there is a reduced decomposition of $\tau$ that contains a reduced decomposition of $\sigma$ as a (not necessarily continuous) substring. 
    Alternatively, the Bruhat order on $S_n$ can be obtained combinatorially by taking the transitive and reflexive closure of the relation: for every $1 \le i < j \le n$, $\s < \s t_{ij} \text{ if } \s(i) < \s(j)$ where $t_{ij}$ is the transposition swapping $i$ and $j$. 
    One can check that this relation indeed is a partial order on $S_n$. 
    We record two basic facts about Schubert varieties of complete flag varieties. 
    \begin{fact}
        Fix a flag $\fflag$ in $\Fl(k^n)$. Let $\s, \tau \in S_n$.
        Then $\s \le \tau$ if and only if $X_{\s, \fflag} \subseteq X_{\tau, \fflag}$.
    \end{fact}
    \begin{fact}
    \label{fct:min-perm}
        Fix a flag $\fflag$ in $\Fl(k^n)$. The dimension of the Schubert cell $X^{\circ}_{\s, \fflag}$ is $\inv(\s)$.
    \end{fact}
    These facts demonstrate that the combinatorics of $S_n$ has tight connections with the geometry of complete flag varieties and Schubert varieties. Indeed, 
    Lakshmibai-Sandhya gave a characterization of singular Schubert varieties and their singular loci using pattern avoidance of indexing permutations: a Schubert variety is singular if and only if the indexing permutation contains patterns $3412$ or $4231$ (in one-line notation) ~\cite{Lakshmibai1990}. Other local properties of Schubert varieties such as being Gorenstein are also characterized using pattern avoidances \cite{Woo2004WHENIA}.
    The main objects that we will study -- the relative Bott-Samelson varieties --  for resolving the singularities of relative Schubert varieties and intersections thereof are also combinatorial in nature.
    
\subsection{Degeneracy loci and relative Schubert varieties}
\label{sec:degeneracy-loci}
    In this section, we recall the notion of the relative position of a pair of flags. These incidence relations involved are determinantal, and thus define schemes that are \textit{degeneracy loci}.
    \begin{proposition}
    \label{prop:rel}
        For $\pflag, \qflag \in \Fl(H)$, there exists a unique permutation in $S_n$ such that there exists an ordered basis $b_1, \ldots, b_n$ of $H$ with the property that
        \[
        b_i \in (P_i \setminus P_{i-1}) \cap (Q_{\s(i)} \setminus Q_{\s(i)-1}).
        \]
    \end{proposition}
    
    \begin{definition}
        For $\pflag, \qflag \in \Fl(H)$,
        we define the \textbf{relative position} of $\pflag$ and $\qflag$ as the unique permutation $\s$ produced in Proposition \ref{prop:rel}, denoted as $r_H(\pflag, \qflag)$, or simply $r(\pflag, \qflag)$ if the vector space is clear from context.
    \end{definition}
    
    \begin{example}
        Let $\omega$ be the order-reversing permutation.
        Two flags $\pflag, \qflag \in \Fl(H)$ are transverse, meaning that \[\dim P_i \cap Q_{n - i} = 0 \text{ for all } i \in [n-1],\] if and only if $r(\pflag, \qflag)$ is $\omega$. 
    \end{example}
    
    \begin{example}
        Two flags $\pflag, \qflag \in \Fl(H)$ are called almost transverse 
        if there exists ${t \in [n-1]}$ such that
        \[
        \dim P_i \cap Q_{n-i} = \begin{cases}
        0 & \text{ if } i \ne t, \\
        1 & \text{ if } i = t.
        \end{cases}
        \] Two flags are almost transverse if and only if $r(\pflag, \qflag) = \omega s_{i}$, where $s_{i}$ is the simple transposition exchanging $i$ and $i + 1$. 
    \end{example}
    
    From the definition, one can easily check that the relative position of two flags is ``anti-symmetric" in the following sense.
    \begin{proposition}
        Let $\pflag, \qflag \in \Fl(H)$.
        Then $r_{H}(\pflag, \qflag) = r_{H}(\qflag, \pflag)^{-1}$.
    \end{proposition}
      
    We now introduce the degeneracy loci in the interest of this paper. 
    \begin{definition}
    Let $\s \in S_n$. Let $S$ a finite-type $k$-scheme carrying a vector bundle $\calH$ of rank $n$. 
    Let $\Fl(\calH)$ be the flag bundle associated with $\calH$ over $S$, with the structure map $\varphi: \Fl(\calH) \to S$.
    Let $p, v$ be sections $S \to \Fl(\calH)$ and denote the fibers of $p, v$ over a point $x \in S$ by $p_x, v_x$. 
    The \textbf{degeneracy locus} $D_{\s}(p; v)$ is a subscheme of $S$
    \[
    D_{\s}(p; v) = \{ x \in S: r(p_x, v_x) \le \s\}.
    \]
    \end{definition}
    Such incidence relations in a family turn out to be determinantal, thereby giving scheme structures on degeneracy loci. 
    We are interested in relative Schubert varieties over some base scheme, which are examples of degeneracy loci, as follows. 
    \begin{example}
    \label{ex:relative-schubert-varieties}
        Let $\s \in S_n$. Let $\varphi^{\ast}\calH$ be the pullback of $\calH$ over $\Fl(\calH)$ and let $\Fl(\varphi^{\ast}\calH)$ be the flag bundle associated with $\varphi^{\ast}\calH$. 
        Let $t: \Fl(\calH) \to \Fl(\varphi^{\ast}\calH)$ be the tautological section. 
        Given a section $p: S \to \Fl(\calH)$, 
        we obtain the section $p': \Fl(\calH) \to \Fl(\varphi^{\ast}\calH)$ canonically induced by $p$; formally, $p' = \id_{\Fl(\calH)} \times (p \circ \varphi)$.
        The relative Schubert variety is a subscheme of $\Fl(\calH)$
        \[
        X_{\s}(p) := D_{\s}(p'; t) \subseteq \Fl(\calH).
        \]
    \end{example}
    
    \begin{example}
    Let $\s \in S_n$. When $S = \Spec k$, $\calH$ is a fixed vector space $H$ of dimension $n$,  and a section $p: S \to \Fl(H)$ yielding a fixed flag $\pflag \in \Fl(H)$, the relative Schubert variety $X_{\s}(p)$ is the Schubert variety $X_{\s, \pflag} \subseteq \Fl(H)$, defined in Section \ref{sec:flag-schubert-varieties}.
    \end{example}

\subsection{Versality}
\label{sec:versality}
    The geometry of intersections of relative Schubert varieties behaves nicely under certain conditions. 
    Morally speaking, when the defining sections of relative Schubert varieties interact with each other minimally in ``most fibers," 
    the local geometry of the intersection is completely determined by the local geometry of each intersecting relative Schubert variety. 
    In ~\cite{chan2019relative}, Chan and Pflueger use the term \textit{versality} for such a condition on the defining flags, which we recall in a slightly different notation as follows. 
    
    Suppose we are given a finite-type $k$-scheme $S$ with a vector bundle $\calH$ of rank $n$. Denote by $\varphi: \Fl(\calH) \to S$ the flag bundle  associated with $\calH$ over $S$,
    and denote by $\psi: \Fr(\calH) \to S$ the frame bundle associated with $\calH$ over $S$. 
    Let $p_1, \ldots, p_m: S \to \Fl(\calH)$ be sections.
    Then we obtain a morphism $\Phi_{p_1, \ldots, p_m}: \Fr(\calH) \to \Fl(k^n)^{m}$ canonically induced by the sections $p_1, \ldots, p_m$; it is defined on the level of closed points by sending $(s, f)$, where $s \in S$ and $f$ is an isomorphism $k^n \to \calH_s$, to the $m$-tuple $(f^{-1}(p_i)_{s})_{i = 1}^{m}$,
    where $(p_i)_s$ is the fiber of $p_i$ over $s$, and for all $i \in [m]$, $f^{-1}(p_i)_s$ denotes the flag 
    \[
        \{0\} \subset (f^{-1}(p_i)_s)_1 \subset \cdots \subset (f^{-1}(p_i)_s)_n = k^n.
    \]
    Certain properties of the map $\Phi_{p_1, \ldots, p_m}$ capture the relative positions of the family of defining flags $\{p_i\}_{i = 1}^{m}$. 
    
    \begin{definition}[{\cite[Definition 3.1]{chan2019relative}}]
    The family of sections $\{p_i\}_{i = 1}^{m}$ is \textbf{versal} if the morphism \[\Phi_{p_1, \ldots, p_m}: \Fr(\calH) \to \Fl(k^n)^{m}\] is smooth. 
    \end{definition}
    
    We collect facts and examples of versal families, proven by Chan-Pflueger using linear-algebraic criterion in {~\cite[Proposition 3.2]{chan2019relative}}. 
    \begin{fact}[{\cite[Lemma 3.4]{chan2019relative}}]
    \label{fact:tautological-versal-are-versal}
        In the situation of Example \ref{ex:relative-schubert-varieties},
        let $t: \Fl(\calH) \to \Fl(\varphi^{\ast} \calH)$ be tautological, let $\{p_i\}_{i = 1}^{m}$ be a versal family of sections $S \to \Fl(\calH)$, 
        and let $\{p'_i\}_{i = 1}^{m}$ be the canonically induced family of sections $\Fl(\calH) \to \Fl(\varphi^{\ast} \calH)$ defined by $p'_i = \id_{\Fl(\calH)} \times (p_i \circ \varphi)$ for all $i \in [m]$. 
        Then the family $\{p'_i\}_{i = 1}^{m} \cup \{t\}$ is versal. 
    \end{fact}
    \begin{example}
    \label{ex:versal-flags-are-transverse}
        When $S = \Spec k$, $H$ is a fixed vector space, and the sections ${p, q: S \to \Fl(H)}$ give a pair of fixed flags in $\pflag, \qflag \in \Fl(H)$, $\Phi_{p,q}$ is smooth if and only if $\pflag, \qflag$ are transverse. 
    \end{example}
    
    Nonetheless, the two sections can be versal, but give non-transverse flags
    in the fiber over a reduced point.
    \begin{example}
        Suppose $S$ is a $1$-dimensional smooth scheme, and sections ${p, q: S \to \Fl(\calH)}$ give two complete flags that are transverse in all but one fiber over a reduced point $x \in S$. Suppose they are almost transverse over $x$, or have relative position $\omega s_{i}$. We have that 
        at $x \in D_{w s_{i}}(p; q)$, $x$ is smooth in $S$ and the codimension of $D_{w s_{i}}$ at $x$ is $\inv(\omega \omega s_{i}) = 1$. Hence by ~\cite[Lemma 3.6]{chan2019relative}, the two flags are versal over $S$. 
        One example to visualize is a family of two distinct points on $\P^1$ moving along a $1$-dimensional base $S$ and the two points meet over a reduced point $x \in S$. The two flags giving the two points are versal over $x$. 
    \end{example}
    
\subsection{Relative Richardson Varieties}
    Our ultimate goal is to construct a Bott-Samelson type resolution of singularities for relative Richardson varieties, which we now define. 
    
    Let $\s, \tau \in S_n$. Let $p, q$ be sections $S \to \Fl(\calH)$, and by Fact \ref{fact:tautological-versal-are-versal}, we obtain the versal family of sections $t, p'$ and $q'$, where $p'$ and $q'$ are the sections $\Fl(\calH) \to \Fl(\varphi^{\ast}\calH)$ canonically induced by $p$ and $q$ respectively. We denote by $X_{\s}(p)$ and $X_{\tau}(q)$ the relative Schubert varieties, defined in Example ~\ref{ex:relative-schubert-varieties}.
    We now define the \textbf{relative Richardson variety} as the intersection of the two relative Schubert varieties defined with respect to these versal sections
    \[
    R_{\s, \tau}(p, q) = X_{\s}(p) \cap X_{\tau}(q).
    \]
    
    The local properties and cohomology of relative Richardson varieties are directly related to those of the intersecting relative Schubert varieties, shown by Chan and Pflueger in \cite{chan2019relative}. They generalized Knutson-Woo-Yong theorem \cite[Theorem 1.1]{KWY2013singularities}  to $\ell > 2$ and to a relative setting. 
    We record Chan-Pflueger's main results. 
    \begin{thm}[{\cite[Theorem 4.1]{chan2019relative}}]
    
    \label{thm:main-1}
        Suppose $P$ is an \'{e}tale-local property preserved under taking product with affine spaces, and 
        $f_{P, m}$ is an $m$-input function such that for finite-type $k$-schemes $X_1, \ldots, X_{m}$ and $x \in \prod_{i = 1}^{m} X_i$, 
        \[P(x, \prod_{i = 1}^{m} X_i) = f_{P, m}(P(\pi_1(x), X_1), \ldots, P(\pi_{m}(x), X_{m})),\] 
        where $\pi_{j}$ is the 
        projection map $\prod_{i = 1}^{m} X_{i} \to X_j$ for all $j \in [m]$.
        Let $t, p_1, \ldots, p_m$ be versal sections $S \to \Fl(\calH)$, let $D_{\s_i}$ denote the degeneracy locus $D_{\s_i}(p_i; t)$ and let $D_{\s_1, \ldots, \s_m}$ denote the intersection of $D_{\s_i}$ for all $i$.  
        Suppose $y \in D_{\s_{1}, \ldots, \s_{m}}$. Then
        \[P(y, D_{\s_1, \ldots, \s_{m}}) = f_{P, m}(P(y, D_{\s_1}), \ldots, P(y, D_{\s_{m}})).
        \]
    \end{thm}
    This theorem can be directly applied to relative Richardson varieties and yields the following.
    \begin{thm}[{\cite[Theorem 5.3]{chan2019relative}}]
        A relative Richardson variety $R_{\s, \tau}(p, q)$ is normal, Cohen-Macaulay and of pure dimension $\inv(\omega \s) + \inv(\omega \tau)$ in $\Fl(\calH)$. 
        Moreover, the smooth locus of $R_{\s, \tau}(p, q)$ is the intersection of the smooth loci of $X_{\s}(p; t)$ and $X_{\tau}(q; t)$.
    \end{thm}

\subsection{Relative Bott-Samelson varieties}
\label{sec:bsv}
        In this section, we recall the definition of relative Bott-Samelson varieties in {\cite[Appendix C]{Fulton1998Schubert}}. They are resolutions of singularities for relative Schubert varieties.
        Let $\s \in S_n$ with length $\ell$ and let $\rho_{\s} = s_{i_1} \cdots s_{i_{\ell}}$ be a reduced decomposition.
        For each $i_j$, let $\Fl(i_j; \calH)$ be the partial flag bundle associated with $\calH$ over $S$ where every point in the fiber over a point $s$ is a partial flag of $\calH_s$ lacking exactly the $i_{j}$-dimensional subspace.
        There is a canonical projection map $\Fl(\calH) \to \Fl(i_j; \calH)$ and we denote by $Z_{s_{i_j}}$ the product $\Fl(\calH) \times_{\Fl(i_j; \calH)} \Fl(\calH)$. The scheme $Z_{s_{i_j}}$ can be set-theoretically described as 
        \[
        \left\{(s, \eflag, \fflag): s \in S, \eflag, \fflag \in \Fl(\calH_s), r(\eflag, \fflag) \le s_{i_j} \right\} \subseteq \Fl(\calH)_S^2.
        \]
        Consider the product 
        \[
        Z_{\rho_{\s}} := Z_{s_{i_1}} \times_{\Fl(\calH)} \cdots \times_{\Fl(\calH)} Z_{s_{i_{\ell}}} \subseteq \Fl(\calH)^{\ell + 1},
        \]
        where the product of each consecutive pair of $Z_{s_{i_j}}$ and $Z_{s_{i_{j+1}}}$ is taken using the second projection of $Z_{s_{i_j}}$ and the first projection of $Z_{s_{i_{j+1}}}$. 
        It can be set-theoretically described as 
        \[
        \left\{(s, \fflag^0, \ldots, \fflag^{\ell}): s \in S, \fflag^{j} \in \Fl(\calH_s), r(\fflag^{j-1}, \fflag^{j}) \le s_{i_{j}} \text{ for all } j \in [\ell]\right\}. 
        \]
        Now let $p$ be a section $S \to \Fl(\calH)$ and denote by $p'$ the section $\Fl(\calH) \to \Fl(\varphi^{\ast}\calH)$ canonically induced by $p$. 
        Recall that the relative Schubert variety $X_{\id}(p) \subseteq \Fl(\calH)$ consists as closed points of $(s, p_s)$ for all $s \in S$. 
        Denote by $\pi_0$ the first projection of $Z_{\rho_{\s}}$ to $\Fl(\calH)$.  

    \begin{definition}
    \label{def:bott-samelson}
        We define the \textbf{relative Bott-Samelson variety} $Z_{\rho_{\s}}(p)$ as
        \[
        Z_{\rho_{\s}}(p) := 
        \pi_{0}^{-1}X_{\id}(p) \subseteq \Fl(\calH)^{\ell + 1}.
        \] 
    \end{definition}
    
    The scheme $Z_{\rho_{\s}}(p)$ can be described set-theoretically as 
        \[
        \left\{(s, \fflag^0, \ldots, \fflag^{\ell}): s \in S, \fflag^{j} \in \Fl(\calH_s), p_s = \fflag^0, r(\fflag^{j-1}, \fflag^{j}) \le s_{i_{j}} \text{ for all } j \in [\ell]\right\}. 
        \]
    
    \begin{proposition}
    \label{prop:rbsv-smooth}
        The relative Bott-Samelson variety $Z_{\rho_{\s}}(p)$ is smooth. 
    \end{proposition}
    
    \begin{proof}
        Since $Z_{\rho_{\s}}(p)$ is a family of  iterated $\P^1$-bundles over $X_{\id}(p)$, $\pi_0$ is a smooth morphism. In addition, $X_{\id}(p)$ is a smooth scheme; therefore, $Z_{\rho_{\s}}(p)$ is smooth.
    \end{proof}
    
    \begin{example}
    \label{ex:absolute-bott-samelson-varieties}
        When $S = \Spec k$ and $H$ is a fixed $n$-dimenisonal vector space over $k$, we recover the absolute Bott-Samelson variety. 
        As an example, let $\s = 4231$ with the reduced decomposition $\rho_{\s} = (34)(12)(23)(34)(12)$. 
        A fixed section $p: S \to \Fl(H)$ gives a fixed flag $\pflag$, and the Bott-Samelson variety $Z_{\rho_{\s}}(p)$ is the subscheme of $\Fl(H)^{6}$ whose construction is schematically depicted in Figure \ref{fig:4231-bsv}. It is a sequence of iterated $\P^1$-fibrations over a point in the flag variety.
        \begin{center}
        \begin{figure}
            \centering \includegraphics[scale=0.4]{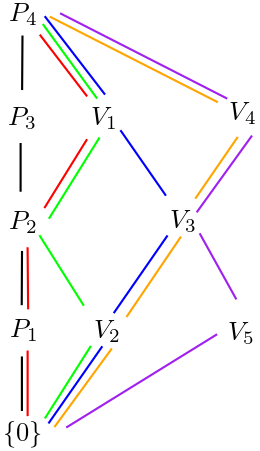}
            \caption{The Bott-Samelson variety $Z_{\rho_{\s}}(p)$ for $\rho_{\s} = (34)(12)(23)(34)(12)$, $\s = 4231$ and a fixed flag $\pflag$. Every point of the variety can be described as a sequence of choices of a red flag, a green flag, a black flag, an orange flag and a purple flag. Adjacent flags differ by at most one subspace $V_j$ of dimension $i_j$ given by $\rho_{\s}$.}
            \label{fig:4231-bsv}
        \end{figure}
        \end{center}
    \end{example}
    In the special case above, $Z_{\rho_{\s}}(p)$ admits a canonical proper morphism to $\Fl(H)$ and is birational to the image -- the Schubert variety $X_{\s, \pflag}$. Therefore, a Bott-Samelson variety is a resolution of singularities for its associated Schubert variety. 
    Bott-Samelson varieties were first introduced by Bott and Samelson ~\cite{bott1958} and later explicitly constructed by Demazure ~\cite{demazure1974} and Hansen ~\cite{hansen1974} as resolutions of singularities for Schubert varieties in complete flag varieties.
    Zelevinsky ~\cite{zelevinskii1983small} generalized the construction to resolve singularities for Grassmannian Schubert varieties. 
    Despite having fibers of large dimenisions and hence not a small resolution, the Bott-Samelson resolutions play a central role in the studies of the geometry of Schubert varieties.
    
    When two Schubert varieties are defined with respect to transverse flags, the product of their Bott-Samelson resolutions over the complete flag variety resolves the singularities of their intersection. Formally, suppose $S = \Spec k$, $H$ is 
        a fixed $n$-dimensional vector space,
        two versal sections $p, q: S \to \Fl(H)$ give transverse flags $\pflag, \qflag \in \Fl(\calH)$, and $\s, \tau \in S_n$ with reduced decompositions $\rho_{\s}, \rho_{\tau}$.
        The variety $Z_{\rho_{\s}}(p) \times_{\Fl(H)} Z_{\rho_{\tau}}(q)$ is smooth and
        the canonical morphism 
        \[
        Z_{\rho_{\s}}(p) \times_{\Fl(H)} Z_{\rho_{\tau}}(q) \to X_{\s, \pflag} \cap X_{\tau, \qflag} \subseteq \Fl(H)
        \] is proper and birational.
        Thus $Z_{\rho_{\s}}(p) \times_{\Fl(H)} Z_{\rho_{\tau}}(q)$ is a resolution of singularities for the Richardson variety $X_{\s, \pflag} \cap X_{\tau, \qflag}$ {\cite[Proof of Theorem 4.2.1]{Brion2005}}.
    We will generalize this result to products of relative Bott-Samelson varieties defined with respect to versal sections.

\section{Products of Relative Bott-Samelson Varieties}
    \subsection{Resolutions for relative complete flag Richardson varieties}
    In this section we prove our main theorem, restated as follows. 
    \begin{thm}
    \label{thm:main}
        Let $S$ be a smooth finite-type $k$-scheme that carries a rank-$n$ vector bundle $\calH$. 
        Let $\s, \tau \in S_n$ with lengths $\ell$ and $\ell'$ respectively. Let $\rho_{\s}$ and $\rho_{\tau}$ be reduced decompositions for $\s$ and $\tau$ respectively.
        Let $p, q: S \to \Fl(\calH)$  be versal sections. Let $X_{\s}(p)$ and $X_{\tau}(q)$ be the relative Schubert varieties defined with respect to $p$ and $q$.
        Let $Z_{\rho_{\s}}(p)$ and $Z_{\rho_{\tau}}(q)$ be the relative Bott-Samelson varieties defined with respect to $p$ and $q$. 
        The fiber product
        \[Z_{\rho_{\s}}(p) \times_{\Fl(\calH)} Z_{\rho_{\tau}}(q)\]
        \begin{enumerate}
            \item is smooth; 
            \item and the projection to the relative Richardson variety \[
            Z_{\rho_{\s}}(p) \times_{\Fl(\calH)} Z_{\rho_{\tau}}(q) \to X_{\s}(p) \cap X_{\tau}(q)
            \] is proper and birational.
        \end{enumerate}
        Therefore, $Z_{\rho_{\s}}(p) \times_{\Fl(\calH)} Z_{\rho_{\tau}}(q)$ is a resolution of singularities of the relative Richardson variety.
    \end{thm}
    
    The crux of the proof is to show that $Z_{\rho_{\s}}(p) \times_{\Fl(\calH)} Z_{\rho_{\tau}}(q)$ is smooth, which we accomplish in two steps. First, we relate the local geometry of $Z_{\rho_{\s}}(p) \times_{\Fl(\calH)} Z_{\rho_{\tau}}(q)$ with that of a principal $\GL_n(k)$-bundle over $Z_{\rho_{\s}}(p) \times_{\Fl(\calH)} Z_{\rho_{\tau}}(q)$.
    Second, we relate the local geometry of the principal bundle with that of relative Bott-Samelson varieties over the complete flag variety, of whose local geometry we have a good grasp, thanks to the versal defining flags. 
    
    We fix notations. 
    Let $\Fl = \Fl(k^n)$, 
    let $\Phi_{p, q}: \Fr(\calH) \to \Fl^2$ be the map canonically induced by $p$ and $q$ in the way described in Section ~\ref{sec:versality}, and denote $Z_{\rho_{\s}}(p) \times_{\Fl(\calH)} Z_{\rho_{\tau}}(q)$ by $Z_{\s, \tau}(p, q)$.
        
    For step one, we construct a principal $\GL_n(k)$-bundle over $Z_{\s, \tau}(p, q)$ as follows. 
    Let $\rho_{\tau}^{-1}$ be the reverse of $\rho_{\tau}$ and $\rho_{\s} \rho_{\tau}^{-1}$ be the concatenation of $\rho_{\s}$ and $\rho_{\tau}^{-1}$. 
    Let $t$ be the tautological section of the trivial bundle $\Fl(\Fl \times k^n) \to \Fl$ over $\Fl$. 
    Consider the scheme $Z_{\rho_{\s}\rho_{\tau}^{-1}}(t) \subseteq \Fl^{\ell + \ell' + 1}$.
    This scheme projects to $\Fl^2$ via the two outer projections, denoted by $\pi_0$ and $\pi_{\ell + \ell'}$ and hence we obtain the fiber product $Z_{\rho_{\s}\rho_{\tau}^{-1}}(t) \times_{\Fl^2} \Fr(\calH)$. There exists a canonical morphism \[\a: Z_{\rho_{\s}\rho_{\tau}^{-1}}(t) \times_{\Fl^2} \Fr(\calH) \to Z_{\s, \tau}(p, q)\] that sends a closed point
    \[
    (\fflag^{\s, 0}, \ldots, \fflag^{\s, \ell} = \fflag^{\tau^{-1}, 0}, \ldots, \fflag^{\tau^{-1}, \ell'}, s, f), 
    \] where $s \in S$, $f: k^n \to \calH_s$ is a vector-space isomorphism and $\fflag^{-,-} \in \Fl(\calH_s)$, to
    \[
    (f(\fflag^{\s, 0}), \ldots, f(\fflag^{\s, \ell})), (f(\fflag^{\tau^{-1}, \ell'}), \ldots, f(\fflag^{\tau^{-1}, 0})).
    \]
    Denote by $\beta$ the canonical proper projection $Z_{\rho_{\s}\rho_{\tau}^{-1}}(t) \times_{\Fl^2} \Fr(\calH) \to \Fr(\calH)$.
    Then we obtain the following diagram
    \begin{center}
        \begin{tikzcd}
            Z_{\s, \tau}(p,q) \ar[d] & Z_{\rho_{\s}\rho_{\tau}^{-1}}(t) \times_{\Fl^2} \Fr(\calH) \ar[d, "\beta"] \ar[l, "\a"] \ar[r] & Z_{\rho_{\s}\rho_{\tau}^{-1}}(t) \ar[d, "{\pi_{0}, \pi_{\ell + \ell'}}"] \\
            S & \Fr(\calH) \ar[l] \ar[r, "\Phi_{p, q}"] & \Fl^2 .
    \end{tikzcd}
    \end{center}
    The intuition behind considering the fiber product $Z_{\rho_{\s}\rho_{\tau}^{-1}}(t) \times_{\Fl^2} \Fr(\calH)$ is that it is the family of sequences of flags where the $(\ell+1)$-th flag has relative position at most $\s$ with $p$ and at most $\tau$ with $q$, upon some change of basis, and  $Z_{\rho_{\s} \rho_{\tau}^{-1}}(t)$ is smooth, by Proposition \ref{prop:rbsv-smooth}. Furthermore, by construction, we have the following lemma. 
        
    \begin{lemma}
        \label{lem:bsv-fiber-product-is-g-torsor}
        The scheme 
        \[
        Z_{\rho_{\s}\rho_{\tau}^{-1}}(t) \times_{\Fl^2} \Fr(\calH)
        \] equipped with the morphisms $\a, \beta$, is the fiber product of $Z_{\s, \tau}(p,q)$ with $\Fr(\calH)$. 
        Hence $Z_{\rho_{\s}\rho_{\tau}^{-1}}(t) \times_{\Fl^2} \Fr(\calH)$ is  isomorphic to $Z_{\s, \tau}(p, q) \times_{S} \Fr(\calH)$, a $\GL_n(k)$-principal bundle over $Z_{\s, \tau}(p, q)$, via a unique isomorphism.
    \end{lemma}
    
    Now we are ready to prove the theorem. 
    \begin{proof}[Proof of Theorem \ref{thm:main}]
        To show smoothness of $Z_{\s, \tau}(p, q)$ in (i), let $x \in Z_{\s, \tau}(p, q)$ and we show that $x$ is a smooth point. 
        By Lemma \ref{lem:bsv-fiber-product-is-g-torsor}, 
        $Z_{\rho_{\s}\rho_{\tau}^{-1}}(t) \times_{\Fl^2} \Fr(\calH)$ is locally isomorphic to $Z_{\s, \tau}(p, q) \times \GL_{n}(k)$. 
        Therefore, if $y \in Z_{\rho_{\s}\rho_{\tau}^{-1}}(t) \times_{\Fl^2} \Fr(\calH)$ such that $\a(y) = x$, 
        then $x$ is smooth in $Z_{\s, \tau}(p, q)$ if and only if $y$ is smooth in $Z_{\rho_{\s}\rho_{\tau}^{-1}}(t) \times_{\Fl^2} \Fr(\calH)$.
        Since $\Phi_{p, q}$ is a smooth morphism, 
        and smoothness of morphisms is stable under base change,
        the morphism $Z_{\rho_{\s}\rho_{\tau}^{-1}}(t) \times_{\Fl^2} \Fr(\calH) \to Z_{\rho_{\s}\rho_{\tau}^{-1}}(t)$ is smooth. 
        Since $Z_{\rho_{\s}\rho_{\tau}^{-1}}(t)$ is a smooth scheme over $k$, 
        and composition of smooth morphisms is smooth, 
        $Z_{\rho_{\s}\rho_{\tau}^{-1}}(t) \times_{\Fl^2} \Fr(\calH)$ is smooth. 
        Therefore, $x$ is a smooth point in $Z_{\s, \tau}(p, q)$ and $Z_{\s, \tau}(p, q)$ is smooth.
        
        For (ii), the scheme $Z_{\s, \tau}(p, q)$ admits a canonical morphism $\pi$ to $\Fl(\calH)$, whose scheme-theoretic image is the relative Richardson variety $X_{\s}(p) \cap X_{\tau}(q)$. 
        Since $Z_{\s, \tau}(p, q)$ is a proper scheme, and $X_{\s}(p) \cap X_{\tau}(q)$ is separated, the morphism $\pi$ is proper. 
        Furthermore, let $R^{\circ}_{\s, \tau}$ be the open subscheme of $X_{\s}(p) \cap X_{\tau}(q)$ where the inequalities of relative positions of flags in all fibers are equalities. 
        It is straightforward to see that $\pi$ is an isomorphism on $R^{\circ}_{\s, \tau}$. Therefore, $\pi$ is birational. 
        
       We conclude that $Z_{\s, \tau}(p, q)$ is a resolution of singularities for the relative Richardson variety $X_{\s}(p) \cap X_{\tau}(q)$. 
    \end{proof}
    
    One immediate corollary is the special case when $S = \Spec k$ and sections $p, q$ give transverse flags $\pflag, \qflag$ of a vector space $H$, previously proven by Brion involving Kleiman's transversality theorem. 
    \begin{corollary}
        Given permutations $\s, \tau \in S_n$ with respective reduced decomposition $\rho_{\s}$, $\rho_{\tau}$ and transverse flags $\pflag, \qflag$ of a vector space $H$, 
        the product $Z_{\rho_{\s}}(\pflag) \times_{\Fl} Z_{\rho_{\tau}}(\qflag)$ is the resolution of singularities for the Richardson variety $X_{\s, \pflag} \cap X_{\tau, \qflag}$.
    \end{corollary}
    
    \subsection{Resolutions for relative Grassmannian Richardson varieties}
    \label{sec:product-grass}
    We have a similar theorem in the case of the versal intersection of two relative Grassmannian Schubert varieties for an application of resolving the singularities of Brill Noether varieties with imposed ramification on twice-marked elliptic curves; see Section \ref{sec:res-bn}. 
    Given integers $n \ge r \ge 1$, recall that $\Gr(r, n)$ is the Grassamannian parametrizing $r$-dimensional subspaces in an $n$-dimensional vector space. A partition $\lambda = (\lambda_1, \ldots, \lambda_r)$ of $\sum \lambda_i$ is \textbf{$\Gr(r, n)$-admissible} if \[
    n - r \ge \lambda_1 \ge \cdots \ge \lambda_r \ge 0. 
    \]
    Alternatively, we can describe $\lambda$ by collecting parts of the same size; that is, we write 
    \[
    \lambda = (\mu_1^{i_1}, \ldots, \mu_j^{i_{j}}) 
    \] where \[
    n - k \ge \mu_1 > \cdots > \mu_{j} > 0
    \] and $\sum i_j = r$. If a partition $\lambda = (\mu_1^{i_1}, \ldots, \mu_{j}^{i_j})$, we say $\lambda$ is of type $j$. 
    For example, the partition $\lambda = (4, 3, 3, 2, 1)$ is a $\Gr(5, 9)$-admissible partition. It can be written as $\lambda = (4^1, 3^2, 2, 1)$ with $(\mu_1, i_1) = (4, 1), (\mu_2, i_2) = (3, 2), (\mu_3, i_3) = (2, 1), (\mu_4, i_4) = (1, 1)$, and is thus of type $4$. 
    \begin{definition}
    \label{def:grass-schubert}
        Given integers $n \ge r \ge 1$, a smooth $k$-scheme $S$ with a vector bundle $\calH$ of rank $n$, $\varphi: \Gr(r, \calH) \to S$ the Grassmannian bundle over $S$, a section ${p: S \to \Fl(\calH)}$ giving a complete flag of subbundles of $\calH$ and a $\Gr(r, n)$-admissible partition $\lambda$,
        the \textbf{relative Grassmannian Schubert variety} \[
        X_{\lambda}(p) := 
        \left\{(x, V \subseteq \calH_{x}): \dim V \cap (p_x)_{n - r + i - \lambda_i} \ge i \text{ for all } 1 \le i \le r \right\} \subseteq \Gr(r, \calH). 
        \]
    \end{definition}
    The incidence conditions are determinantal, thus yielding subschemes of $\Gr(r, \calH)$. 
    
    \begin{definition}
    \label{def:res-gr}
        Given $n \ge r \ge 1$, a smooth $k$-scheme $S$ with a vector bundle $\calH$ of rank $n$, a section $p: S \to \Fl(\calH)$ giving a complete flag of subbundles of $\calH$, a $\Gr(r, n)$-admissible partition $\lambda$ written as $(\mu_1^{i_1}, \ldots, \mu_j^{i_j})$, and $a_s = \sum_{\ell = 1}^{s} i_\ell$,
        the \textbf{relative Bott-Samelson resolution for relative Grassmannian Schubert varieties} is
        \[
        Z_{\lambda}(p) = \left\{(x, V_1 \subset \cdots \subset V_j): \dim V_s = a_s, V_s \subseteq (p_x)_{n-r +a_s - \lambda_{a_s}} \text{ for all } 1 \le s \le j \right\} \subseteq \Fl(a_1, \ldots, a_j; \calH).
        \]
    \end{definition}
    The incidence relations are again determinantal, thus yielding a variety. This description is a generalization of the Bott-Samelson/Zelevinsky resolution to the relative setting; see \cite{Coskun2011RigidAN} for a definition in the absolute case. 

    \begin{thm}
    \label{thm:main-gr}
        Given integers $n \ge r \ge 1$, a smooth $k$-scheme $S$ with a vector bundle $\calH$ of rank $n$, the flag bundle $\Fl(\calH) \to S$ associated with $\calH$, the Grassmannian bundle $\Gr(r, \calH) \to S$, two versal sections $p, q: S \to \Fl(\calH)$ and two $\Gr(r, n)$-admissible partitions $\lambda, \lambda'$,
        the product of relative Bott-Samelson resolutions for relative Grassmannian Schubert varieties over the Grassmann bundle 
        \[
        Z_{\lambda}(p) \times_{\Gr(r, \calH)} Z_{\lambda'}(q) 
        \] 
        \begin{enumerate}
            \item is smooth; and 
            \item the projection to the relative Grassmannian Richardson variety
            \[
            Z_{\lambda}(p) \times_{\Gr(r, \calH)} Z_{\lambda'}(q)  \to X_{\lambda}(p) \cap X_{\lambda'}(q) \subseteq \Gr(r, \calH)
            \] is proper and birational. 
        \end{enumerate}
        Therefore, $Z_{\lambda}(p) \times_{\Gr(r, \calH)} Z_{\lambda'}(q)$ is a resolution of singularities for the relative Grassmannian Richardson variety. 
    \end{thm}
    
    \begin{proof}
        We exploit the idea in the proof of Theorem \ref{thm:main} and 
        denote ${Z_{\lambda}(p) \times_{\Gr(r, \calH)} Z_{\lambda'}(q)}$ by $Z_{\lambda, \lambda'}(p, q)$.
        To show smoothness of $Z_{\lambda, \lambda'}(p, q)$, we construct a principal $\GL$-bundle over $Z_{\lambda, \lambda'}(p, q)$.
        Suppose $\lambda$ and $\lambda'$ are of type $j$ and $j'$ respectively, and define $Z_{\lambda, \lambda'}$ as follows: 
        \begin{align*}
            Z_{\lambda, \lambda'} = \left\{F_1, \ldots, F_n, V_1, \ldots, V_j, U_1, \ldots, U_{j'},
            G_1, \ldots, G_n: (\ast)\right\}
        \end{align*}
        where the condition $(\ast)$ consists of the followings:
        \begin{enumerate}
            \item $F_1 \subset F_2 \subset \cdots \subset F_{n-1} \subset F_n$ and $G_1 \subset G_2 \subset \cdots \subset G_{n-1} \subset G_n$ are complete flags in $\Fl$; 
            \item $V_1 \subset V_2 \subset \cdots \subset V_{j-1} \subset V_j$ and $U_1 \subset U_2 \subset\cdots \subset U_{j'-1} \subset U_{j'}$; 
            \item for all $1 \le \ell \le j$, $\dim V_{\ell} = a_{\ell}$ and $V_{\ell} \subseteq F_{n - r + a_{\ell} - \lambda_{a_{\ell}}}$;
            \item for all $1 \le \ell' \le j'$, $\dim U_{\ell'} = b_{\ell}$ and $U_{\ell'} \subseteq G_{n - r + b_{\ell'} - \lambda'_{a_{\ell'}}}$; and 
            \item $V_j = U_{j'}$. 
        \end{enumerate}
        The incidence relations are determinantal, thus defining a variety. Furthermore, $Z_{\lambda, \lambda'}$ is smooth because it is an iterated tower of Grassmannian bundles. It admits projections $\pi, \pi'$ to $\Fl^2$ by projecting the complete flags $\{0\} \subset F_1 \subset \cdots \subset F_n = k^n$ and ${\{0\} \subset G_1 \subset \cdots \subset G_n = k^n}$. 
        Again by construction, the product $Z_{\lambda, \lambda'} \times_{\Fl^2} \Fr(\calH)$ is isomorphic via a unique isomorphism to the product $Z_{\lambda, \lambda'}(p, q) \times_{S} \Fr(\calH)$, and  the following diagram commutes:
        \begin{center}
        \begin{tikzcd}
            Z_{\lambda, \lambda'}(p, q) \ar[d] & Z_{\lambda, \lambda'} \times_{\Fl^2} \Fr(\calH) \ar[d, "\beta"] \ar[l, "\a"] \ar[r] & Z_{\lambda, \lambda'} \ar[d, "{\pi, \pi'}"] \\
            S & \Fr(\calH) \ar[l] \ar[r, "\Phi_{p, q}"] & \Fl^2 .
        \end{tikzcd}
        \end{center}
        Hence $Z_{\lambda, \lambda'} \times_{\Fl^2} \Fr(\calH)$ is a principal $\GL$-bundle over $Z_{\lambda, \lambda'}(p, q)$.
        As similarly argued in the proof of Theorem \ref{thm:main},
        $Z_{\lambda, \lambda'}(p, q)$ is smooth if and only if $Z_{\lambda, \lambda'}$ is smooth, which holds true as previously discussed. 
        
        The variety $Z_{\lambda, \lambda'}(p, q)$ admits a canonical projection $\pi_r$ to $\Gr(r, \calH)$ whose scheme-theoretic image is the relative Grassmannian Richardson variety $X_{\lambda}(p) \cap X_{\lambda'}(q)$. Since $Z_{\lambda, \lambda'}(p, q)$ is proper and $X_{\lambda}(p) \cap X_{\lambda'}(q)$ is separated, the morphism $\pi_r$ is proper. 
        Let $R^{\circ}_{\lambda, \lambda'}$ be the open subscheme of $X_{\lambda}(p) \cap X_{\lambda'}(q)$ where the inequalities of relative positions are equalities. Since $\pi_r$ is an isomorphism on $R^{\circ}_{\lambda, \lambda'}$, $\pi_r$ is birational. 
     \end{proof}
     
    We give the essentially smallest interesting example in the case of relative Grassmannian Schubert intersections. 
    \begin{example}
    \label{ex:grass}
        Let $r = 2, n = 4$. 
        Let $S$ be a $1$-dimensional smooth $k$-scheme with a rank-$4$ vector bundle $\calH$. Further suppose we are given versal sections $p, q$ such that over a unique reduced point $x \in S$, 
        \[
        \dim (p_{x})_{4 - j} \cap (q_{x})_j = \begin{cases}
        0, & \text{if }j = 1, 3, 4,\\
        1, & \text{if }j = 2,
        \end{cases}
        \]
        and that $p_s$ and $q_s$ are transverse everywhere else; see Figure \ref{fig:versal}.
        Let $\lambda = \lambda'$ be the $\Gr(2, 4)$-admissible partition $(1, 0)$ and for any $s \in S$, let $P_s$ and $Q_s$ denote the projectivization of $(p_{s})_2$ and $(q_s)_2$ respectively. 
        Then the relative Grassmannian Richardson variety ${X_{\lambda}(p) \cap X_{\lambda'}(q)}$ parametrizes $(s, L)$ where $s \in S$, $L$ is a line in $\P(\calH_{s})$ and $L$ intersects the two lines $P_s$ and $Q_s$ each at a point.
        \begin{figure}[h]
            \centering
            \includegraphics[scale=0.5]{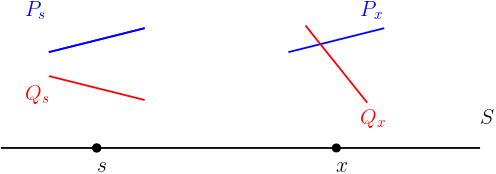}
            \caption{A picture depicting the versality of $p$ and $q$ over a $1$-dimensional scheme $S$: in $\P(\calH_s)$ for $s \ne x$, $P_s$, $Q_s$ are two skew lines in $\P(\calH_s)$; and in $\P(\calH_x)$, $P_x$ and $Q_x$ are two lines touching at a point.}
            \label{fig:versal}
        \end{figure}
        By \cite{chan2019relative}, the fiber of $X_{\lambda}(p) \cap X_{\lambda'}(q)$ over $x$ is the union of two $\P^2$ intersecting at $\P^1$, and the singular locus of $X^r_{\s}(p) \cap X^{r}_{\tau}(q)$ is precisely the the points $(x, L = P_x)$ and $(x, L = Q_x)$ contained in the $\P^1$ intersection. The smooth points of $X_{\lambda}(p) \cap X_{\lambda'}(q)$ are $(s, L)$ in one of the following situations:
            \begin{enumerate}
            \item $s = x$, $L \cap P_x = L \cap Q_x = P_x \cap Q_x$; 
            \item $s = x$, $L \cap P_x \ne P_x \cap Q_x$ and $L \cap Q_x \ne P_x \cap Q_x$;  
            \item $s \ne x$.
            \end{enumerate}
            See Figure \ref{fig:sing} and Figure \ref{fig:smooth} for pictures picturing the incidence relations of $L$ and $P_s, Q_s$ depending on the smoothness at $(s, L)$ in $X^r_{\s}(p) \cap X^{r}_{\tau}(q)$.
        \begin{figure}[h]
            \centering
            \includegraphics[scale=0.3]{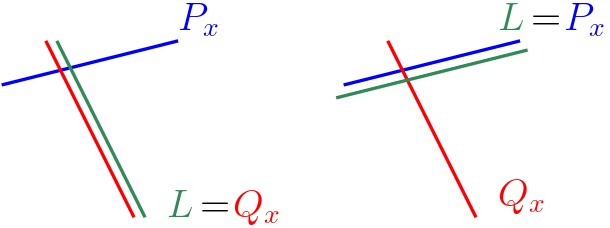}
            \caption{A picture illustrating the incidence relations of $L$ and $P_x, Q_x$ in $\P(\calH_x)$ when $(x, L)$ is a singular point in $X_{\lambda}(p) \cap X_{\lambda'}^{r}(q)$; that is, precisely when $L = P_x$ or $L = Q_x$.}
            \label{fig:sing}
        \end{figure}
        
        \begin{figure}[h]
            \centering \includegraphics[scale=0.4]{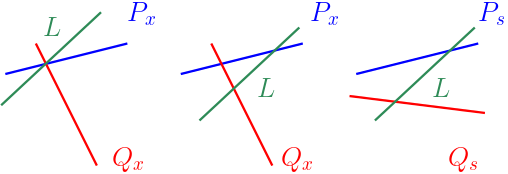}
            \caption{A picture illustrating the incidence relations of $L$ and $P_s, Q_s$ in $\P(\calH_s)$ when $(s, L)$ is a smooth point in $X_{\lambda}(p) \cap X_{\lambda'}(q)$. }
            \label{fig:smooth}
        \end{figure}
        
        The smooth variety $Z_{\lambda}(p) \times_{\Gr(r, \calH)} Z_{\lambda'}(q)$ parametrizes $(s, L, p_1, p_2, H_1, H_2)$ where the points $p_1, p_2$ are contained in $L \cap P_x, L \cap Q_x$ respectively, and the planes $H_1, H_2$ contain $\{L, P_s\}$, $\{L, Q_s\}$ respectively. 
        Over $(s \ne x, L)$, $(x, L \ne P_x)$ or $(x, L \ne Q_x)$, the line $L$ uniquely determines $p_1, p_2, H_1$ and $H_2$ because $p_1 = L \cap P_x$, $p_2 = L \cap Q_x$,  and $H_1, H_2$ are the spans of $\{L, P_x\}$ and $\{L, Q_x\}$ respectively; therefore, the projection $\pi$ is an isomorphism on the open set of $X_{\lambda}(p) \cap X_{\lambda'}(q)$ where equalities of relative positions hold; see Figure \ref{fig:bsv-smooth}.
        \begin{figure}[h]
            \centering
            \includegraphics[scale=0.5]{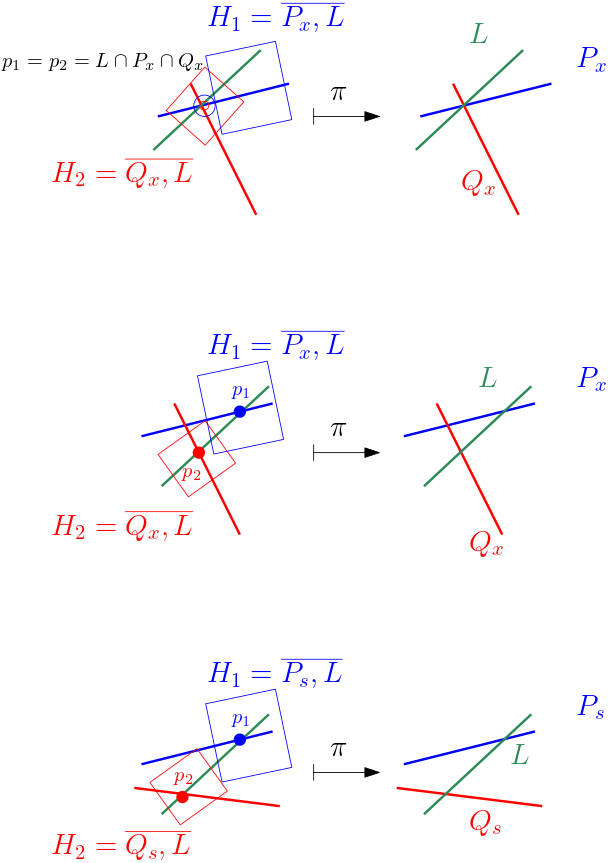}
            \caption{A picture showing $(L, p_1, p_2, H_1, H_2)$ in $\P(\calH_s)$ for  smooth $(s, L)$ in $X_{\lambda}(p) \cap X_{\lambda'}(q)$. The points $p_1, p_2$ are indicated by the blue point and the red point respectively; the planes $H_1, H_2$ are indicated by the blue and red planes respectively. Here, $p_1, p_2, H_1, H_2$ are uniquely determined by $L, P_s, Q_s$ in all three situations. The morphism $\pi$ is an isomorphism by forgetting the planes and the points.}
            \label{fig:bsv-smooth}
        \end{figure} Over $(x, P_x)$ or $(x, Q_x)$, the fiber of $\pi$ is  isomorphic to $\P^1 \times \P^1$; see Figure \ref{fig:bsv-sing}. 
        \begin{figure}[h]
            \centering
            \includegraphics[scale=0.3]{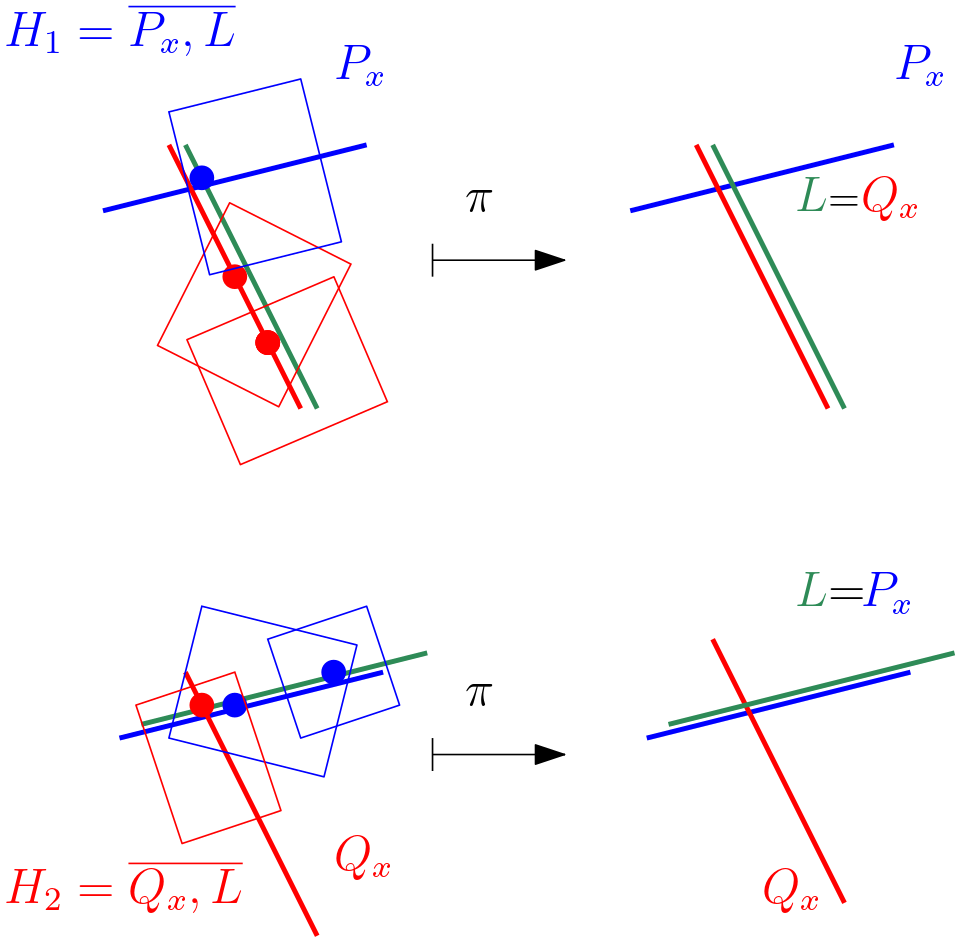}
            \caption{A picture illustrating $(L, p_1, p_2, H_1, H_2)$ in $\P(\calH_x)$ where $(s, L)$ is singular in $X_{\lambda}(p) \cap X_{\lambda'}(q)$. In the top figure when $L = Q_x$, the blue point $p_1 = L \cap P_x$ and blue plane $H_1 = \overline{P_x, L}$; hence they are uniquely determined. However, there are $\P^1$-worth of choices for the red point $p_2$ and $\P^1$-worth choices for the red plane $H_2$. Thus $\pi$ has fiber isomorphic to $\P^1 \times \P^1$. Similarly for the bottom figure for the case when $L = P_x$. }
            \label{fig:bsv-sing}
        \end{figure}
    \end{example}

\section{Application to Brill-Noether Theory}

In this section, we return to Brill-Noether varieties with imposed ramifications, which motivated our consideration of resolution of singularities for relative Richardson varieties. 

\subsection{Resolution of singularities for twice-pointed Brill-Noether varieties}
\label{sec:res-bn}
Let ${g \ge 0}, d > r \ge 0$, and let $(X; P, Q)$ be twice-marked smooth projective curve. Denote by $\vsa$ and $\vsb$ sequences
$0 \le a_1 < \cdots < a_r < a_{r+1}\le d$ and $0 \le b_1 < \cdots < b_r < b_{r+1} \le d$.
The moduli space of linear series of projective rank $r$ of degree $d$ with imposed ramification at $P, Q$ prescribed by $\vsa, \vsb$ respectively is the subscheme of the classical Brill-Noether variety $G^{r}_{d}(X)$
        \[
        G^{r}_{d}(X, (P, \vsa), (Q, \vsb)) \subseteq G^r_d(X);
        \] see \cite{cop}. Chan, Osserman and Pflueger showed that in the case when $g = 1$ and $P-Q$ is not a torsion point of order less than or equal to $d$ in $\Pic^{0}(X)$, the variety $G^{r}_{d}(X, (P, \vsa), (Q, \vsb))$ is the intersection of two relative Grassmannian Schubert varieties $G^r_d(X, (P, \vsa))$ and $G^r_d(X, (Q, \vsb))$ over the base scheme $\Pic^d(X)$ (also see \cite[Corollary 6.2]{chan2019relative}) and the variety has singular locus precisely the union
        \[
        \left(G^{r}_{d}(X, (P, \vsa))^{\text{sing}} \cap G^{r}_{d}(X, (Q, \vsb))\right) \cup \left(G^{r}_{d}(X, (P, \vsa)) \cap G^{r}_{d}(X, (Q, \vsb))^{\text{sing}}\right).
    \]
    Our Theorem \ref{thm:main-gr} can be applied as follows. 
    
    Let $E$ be an elliptic curve marked by $P, Q$ where $P - Q$ is not a torsion point of order weakly less than $d$ in $\Pic^{0}(E)$.
    Let $\calL$ be the Poincar\'{e} line bundle on $E \times \Pic^d(E)$ and let $\pi$ be the projection $E \times \Pic^{d}(E) \to \Pic^{d}(E)$. Then the base scheme $\Pic^{d}(E)$ carries a rank-$d$ vector bundle $\calH = \pi_{\ast}\calL$.
    Over each $[L] \in \Pic^{d}(E)$, set $p_{[L]}, q_{[L]}$ to be the flags of global sections $H^{0}(E, L)$
    \begin{align*}
        p_{[L]} &= \{0\} \subset H^{0}(E, L((1-d)P)) \subset \cdots \subset H^{0}(E, L(-P)) \subset H^{0}(E, L),\\
        q_{[L]} &= \{0\} \subset H^{0}(E, L((1-d)Q)) \subset \cdots \subset H^{0}(E, L(-Q)) \subset H^{0}(E, L).
    \end{align*}
    furnishing two sections $p, q: \Pic^{d}(E) \to \Fl(\calH)$. Let $\vsa$ and $\vsb$ be increasing sequences of $r + 1$ nonnegative numbers no greater than $d$. The sections $p, q$ are proven to be versal in \cite[Lemma 6.1]{chan2019relative}.
    Define $\Gr(r + 1, d)$-admissible partitions $\lambda = (\lambda_i)$ and $\lambda' = (\lambda'_{i'})$  where \begin{align*}\lambda_i &= a_{r + 1 - (i-1)} - (r + 1 - i),\\\lambda'_{i'} &= b_{r + 1 - (i'-1)} - (r + 1 - i')
    \end{align*} for $i, i' \in [r+1]$.
    \begin{thm}
    \label{thm:brill-noether-resing}
        In the situation above, 
        the product of Bott-Samelson resolutions over the Brill-Noether variety $G^{r}_d(E)$
        \[
        Z_{\lambda}(p) \times_{G^{r}_{d}(E)} Z_{\lambda'}(q)
        \] is a resolution of singularities of $G^{r}_{d}(E, (P, \vsa), (Q, \vsb))$. 
    \end{thm}
    \begin{proof}
        This is a direct appliaction of Theorem \ref{thm:main-gr}. 
    \end{proof}
    
    \begin{example}
        In the situation of Theorem \ref{thm:brill-noether-resing}, $g=1$, $d = 4$, $r = 1$, ${\vsa = \vsb = (0, 2)}$ and $\lambda = \lambda' = (1, 0)$, the resolution of singularities $Z_{\lambda}(p) \times_{G^{1}_{4}(E)} Z_{\lambda'}(q)$ for $G^{1}_{4}(X, (P, (0, 2)), (Q, (0, 2)))$ is essentially described as in Example \ref{ex:grass}.
    \end{example}

    \subsection{Open problems}
    \label{sec:open}
    We now state two conjectures which generalize Theorem \ref{thm:main-brill-noether-resing} to Brill-Noether varieties on twice-marked curves of higher genus.
    First we generalize versality to partial flags. 
    Let $n, m \ge 1$, and let $S$ be a finite-type $k$-scheme with a vector bundle $\calH$ of rank $n$. 
    Let $\{d_{\bl}^{i}\}_{i = 1}^{m}$ be increasing sequences of positive integers less than $n$ and let $\Fl(d^i_{\bl}; \calH) \to S$ be the partial flag bundle associated with $\calH$ with dimensions specified by $d^i_{\bl}$ for all $i \in [m]$.
    Let $\Fr(\calH) \to S$ be the frame bundle over $S$ associated with $\calH$. 
    For each $i \in [m]$, suppose $p_i: S \to \Fl(d^{i}_{\bl}: \calH)$ is a section. 
    The family of partial flag sections $p_1, \ldots, p_m$ is \textbf{versal} if the induced morphism 
    \[
    \Phi_{p_1, \ldots, p_m}: \Fr(\calH) \to \prod_{i=1}^{m} \Fl(d^i_{\bl}; \calH)
    \] is smooth. 
    
    Let $n, r \ge 0$, let $\l$ be $\Gr(r, n)$-admissible partitions that can be written as $(\mu_1^{i_1}, \ldots, \mu_{j}^{i_j})$ and set $a_s = \sum_{\ell = 1}^{s} i_{\ell}$. A section $p: S \to \Fl(d_{\bl}; \calH)$ for some increasing sequence of positive integers $d_{\bl}$ less than $n$ is \textbf{$\lambda$-compatible} if
    \[
    (p_s)_{n - r + a_s - \lambda_{a_s}} = d_j \in d_{\bl}
    \] for some $j$ for all $x \in S$. 
    Then define relative Grassmannian Schubert varieties $X_{\lambda}(p)$ with respect to partial flags and Bott-Samelson resolutions $Z_{\lambda}(p)$ for relative Grassmannian Schubert varieties with respect to partial flags similarly to how those are defined in Definition \ref{def:grass-schubert} and \ref{def:res-gr}. Given $\Gr(r, n)$-admissible partitions $\lambda, \lambda'$ and versal sections $p, q$ that are $\lambda$-compatible and $\lambda'$-compatible respectively, we have the following conjecture. 
    \begin{conjecture}
        The product of Bott-Samelson resolutions for relative Grassmannian Schubert varieties with respect to partial flags over $\Gr(r, \calH)$ 
        \[
        Z_{\lambda}(p) \times_{\Gr(r, \calH)} Z_{\lambda'}(q)
        \] is smooth and admits a proper birational morphism to the relative Richardson variety with respect to partial flags $X_{\lambda}(p) \cap X_{\lambda'}(q)$. 
    \end{conjecture}
    
    This conjecture may be specialized to a resolution of singularities for the Brill-Noether varieties on twice-marked curves in higher genus. Let $g > 1$, $d> r \ge 0$ and let $(X, P, Q)$ be a smooth projective curve of genus $g$. For nonnegative integers $N, a, b$ such that ${N, N - a, N-b\ge 2g - 1}$ and $d = N - a - b$, the variety $\Pic^{N}(X)$ carries a vector bundle $\calH$ such that $\calH|_{[L]} = L(aP + bQ)$. 
    Then the Brill-Noether variety $G^{r}_d(X)$ and $G^{r}_{d}(X, (P, \vsa), (Q, \vsb))$ are conjecturally isomorphic to relative Richardson varieties defined with respect to partial flags; see {\cite[Conjecture 6.3]{chan2019relative}}. 
    Suppose $\lambda, \lambda'$ are $\Gr(r+1, d - g + 1)$-admissble partitions given by vanishing sequences $\vsa$ and $\vsb$ as in Section \ref{sec:res-bn}. Then we have a generalization of Theorem \ref{thm:brill-noether-resing}.
    \begin{conjecture}
        In the situation above where $g > 1$, the product of Bott-Samelson resolutions over the Brill-Noether variety $G^{r}_d(X)$ 
        \[
        Z_{\lambda}(p) \times_{G^r_d(X)} Z_{\lambda'}(q)
        \] is a resolution of singularities of $G^{r}_{d}(X, (P, \vsa), (Q, \vsb))$. 
    \end{conjecture}

\bibliographystyle{alpha}
\bibliography{rbsv.bib}

\end{document}